\documentclass[11pt,twoside]{article}
\usepackage{latexsym}
\usepackage{amssymb,amsbsy,amsmath,amsfonts,amssymb,amscd} 
\usepackage{graphicx}
\usepackage{hyperref}
\newcommand{\ind}[1]{\mathbb{I}_{\{#1\} }} 
\newcommand{\fer}[1]{(\ref{#1})}
\newcommand {\wt}[1] {{\widetilde #1}}
\newcommand{\commentout}[1]{}
\newcommand{\RR}{\mathbb{R}}

\newcommand {\al} {\alpha}

\newcommand {\sg} {\sigma}

\newcommand {\lb} {\lambda}
\newcommand {\f} {\frac}
\newcommand {\p} {\partial}
\newcommand {\proof} {\noindent {\bf Proof}. }
\newcommand{\beq}{\begin{equation}}
\newcommand{\eeq}{\end{equation}}

\newtheorem{thm}{Theorem}

\newtheorem{prop}{Proposition}
\newtheorem{lemme}{Lemma}

\newcommand{\qed}{{ \hfill
{\unskip\kern 6pt\penalty 500
\raise -2pt\hbox{\vrule\vbox to 6pt{\hrule width 6pt \vfill\hrule}\vrule} \par} }}
\title{ {\bf Relaxation and self-sustained oscillations in the time elapsed neuron network model} 
\\
}

\author{
Khashayar Pakdaman\thanks{Univ Paris Diderot. Institut Jacques Monod  Email: pakdaman@ijm.univ-paris-diderot.fr}
\and 
Beno\^ \i t Perthame\thanks{
UPMC, CNRS UMR 7598, Laboratoire Jacques-Louis Lions, F-75252, Paris Cedex 05, INRIA-Paris-Rocquencourt EPI BANG and Institut Universitaire de France. Email: benoit.perthame@upmc.fr}
\and
Delphine Salort\footnotemark[1] \thanks{Email: salort.delphine@ijm.univ-paris-diderot.fr}
}

\date{\today}

\begin{document}
\maketitle
\pagestyle{plain}
\pagenumbering{arabic}

\begin{abstract}
The {\em time elapsed model}  describes the firing activity of an homogenous assembly of neurons thanks to the distribution of times elapsed since the last discharge. It gives a mathematical description of the probability density of neurons structured by this time. In an earlier work, based on generalized relative entropy methods, it is proved that for highly or weakly connected networks the model exhibits relaxation to the steady state and for moderately connected networks it is obtained numerical evidence of appearance of self-sustained periodic solutions.

Here, we go further and, using the particular form of the model, we quantify the regime where relaxation to a stationary state occurs in terms of the network connectivity. To introduce our methodology, we first consider the case where the neurons are not connected  and we give a new statement showing that total asynchronous firing of neurons appears asymptotically. In a second step, we consider the case with connections and give a low connectivity condition that  still leads to asynchronous firing. Our low connectivity condition is somehow sharp because we can give an example, when this condition is not fulfilled, where synchronous rhythmic activity occurs. Indeed, we are able to build several explicit  families of periodic solutions. Our construction  is fully nonlinear and the resynchronization of the neural activity in the network does not follow from bifurcation analysis. It relies on an algebraically nonlinear boundary condition that occurs in the model.

These analytic results are compared with numerical simulations under broader hypotheses and shown to be robust.

\end{abstract}

\noindent {\bf Key-words} Neuron networks, desynchronization, self-sustained oscillations

\section{Introduction} 
\label{sec:model}

In nervous systems, neuronal circuits carry out tasks of information transmission and processing. During such functions, neuronal activity is evoked by incoming signals. However, in many assemblies, neuronal electrical activity is present and persists even in the absence of external signals. This form of activity that is not elicited by inputs is referred to as spontaneous activity (SA). That SA plays a pivotal role in functions such as respiratory rhythmogenesis that are and need to be maintained throughout the lifetime of organisms has been long recognized (see for instance \cite{BDC}). Furthermore, for over the past twenty years, experimental evidence has been gathering for the widespread occurrence of SA in different areas of nervous systems, such as the spinal cord, the retina, the cortex and the hippocampus, specially during various phases of development whereby SA mediates neuronal and network maturation (for reviews see \cite{BF}\cite{MB}). In many instances, SA takes
  the form of recurring synchronous activity that emanates from a specific region and subpopulation of neurons and propagates throughout neuronal populations interacting through excitatory connections. In this work, our focus is on the mechanisms through which synchronous neuronal  activity is initiated in tightly connected excitatory neuronal assemblies. 

Building on our previous mathematical analysis of dynamics of neuronal assemblies \cite{Pak_Per_Sal}, our main purpose is to use a mathematically tractable model of large neuronal populations incorporating a minimal set of neuronal properties and construct explicitly periodic solutions of this model corresponding to synchronous oscillations reminiscent of SA. In the process, we also improve previous conditions ensuring asynchrony of assemblies.
Our mathematical analysis and explicit construction are then compared with numerical simulations under broader hypotheses and shown to be robust. Our approach highlights basic mechanisms lying at the core of synchronous rhythm generation in assemblies of interacting noisy excitable systems.

That coupled oscillators can produce collective rhythms is well documented. Indeed, in such systems, isolated units display periodic oscillations. Interactions between such units can entrain them into phase locking yielding global synchronous rhythms at the population level. The situation for coupled excitable systems is different. In this case, isolated units stabilize at a resting state and do not have intrinsic periodicity. It is in the presence of noise that they display occasional irregular discharges. Interactions between units can transform this irregular spontaneous activity of isolated excitable systems into a global coherent and regular rhythm. What makes this form of spontaneous regular synchronous activity of particular interest is that it emerges in a population of excitable rather than oscillating units. In other words, unlike synchrony in coupled oscillators, here synchrony and regularity of the global population are not direct consequences of intrinsic periodicity of the units.

The emergence of such spontaneous synchronous oscillations in populations of noisy excitable neuronal models has been widely observed in the literature. One of the seldom cited earliest reports of this phenomenon is MacGregor and Palasek's study of randomly connected populations of neuromimes \cite{cf:macgregor}. Other examples include the description of similar phenomena in models such as the active rotator \cite{cf:shinomoto1986a}\cite{cf:poquet}, elementary neuron models \cite{cf:pham98b}, the FitzHugh-Nagumo \cite{cf:toral}, the Morris-Lecar \cite{cf:han}, the Hindmarsch-Rose \cite{cf:du}, the Hodgkin-Huxley \cite{cf:wang} and even models implementing detailed biophysical properties \cite{cf:kosmidis}. In some of these references, the units are diffusively coupled, in others they interact through excitatory pulses; in some connectivity is all-to-all, while others deal with random networks. Our enumeration, which does not intend to be exhaustive, illustrates the ease with which  assemblies of excitable units generate noise-induced synchronous regular activity, irrespective of model and network details.

It is the widespread occurrence of the phenomenon in diverse network models that motivated us to design an elementary model with minimal assumptions shared by all excitable systems to gain better understanding of the way such rhythmic patterns arise. Previous explorations of the phenomenon are mainly numerical. Here we present the first example where there is an explicit construction of solutions for this problem. While our model was originally adjusted to experimental data \cite{cf:pham98a}, in its present version, it contains only a bare skeleton of physiological neuronal properties. This reduction serves to emphasize and extract the key mechanisms at work. It shows that even with minimal assumptions it is possible to generate synchronous activity. Another advantage of the explicit construction is that it opens the door to study other forms of patterns that may arise in assemblies of interacting systems. We highlight this aspect by constructing other forms of periodic synchronous oscillations that had not been reported previously. Using numerical simulations, we verify that the results are robust to the details of the model.

From a mathematical standpoint, the time elapsed neuron population model belongs to the class of structured models that appear in a broad range of applications (see the reviews \cite{MD, Pe}). The version considered in this paper possesses a threshold-like firing rate. This is similar to threshold like maturation and death rates in models of cell maturation, fish  populations etc \cite{hbid}. In the same way as for these other models, this leads to a natural relation between our model and equations with state dependent delays. We take advantage of this property in the analysis of the dynamics of populations of unconnected neurons. However, for strongly interacting neuronal populations, no such relations exist and the analysis of the dynamics cannot be performed using this technique. In this respect, the periodic oscillations that are obtained in this work are also novel to structured populations and should be of broader interest in the field of structured models.

This article is organized as follows. In section \ref{sec_model}, we present the time elapsed neuron population model and discuss the underlying assumptions with respect to other related models. In the following section \ref{sec:linear}, we consider the case of unconnected neuronal populations, henceforth referred to as neuronal ensemble. We prove that such a population tends to a state of total asynchronous neuronal activity. In our analysis , we use the underlying delay equation. This section also introduces tools that are used later in the paper. In section \ref{sec:desync}, we extend the conditions for asynchronous activity to interconnected populations of neurons satisfying a criterion on postdischarge recovery rate and the maximal refractory period. After these results, we focus on situations where the network can sustain synchrony. More precisely, in section \ref{sec:periodic}, we build an explicit example of a network that possesses periodic solutions corresponding to rhythmic and synchronous activity in the network. Finally, in section \ref{sec:num}, numerical results and variants of the model show the robustness of the results. We conclude the paper in section \ref{sec:conclusion} by discussing our results, possible extensions and limitations.

\section{The  time elapsed neuron population model}
\label{sec_model}

The population model analyzed in this work encompasses a minimal set of properties reproducing the following aspects of neuronal behavior. Many neurons generate trains of stereotyped electrical pulses --- referred to as spikes or action potentials --- in response to incoming stimulations. Following each discharge, the neuron undergoes a period of refractoriness during which it is less responsive to inputs, before recovering its excitability \cite{cf:kandel}. 

One key aspect of neuronal coding is that the shape of the spike varies little. The main carrier of information is the discharge times or some statistics of the discharge times. Motivated by this, many neuronal and network models neglect the mechanisms underlying spike generation and instead represent neuronal dynamics in terms of discharge times. A general model of neuronal population dynamics based upon occurrence times of events was introduced in the pioneering work of Perkel \cite{perkel76a,perkel76b,perkel76c}. Other related variants are the spike response model or the  integrate and fire models \cite{GK}. In this work, we consider a specific neuronal model that is inspired from our previous studies that aimed to reproduce SA in experimental preparations \cite{cf:pham98a}. 

The fundamental assumption in our model is to describe the postdischarge recovery of the neuronal membranes through an instantaneous firing rate that depends on the time elapsed since the last discharge and the inputs by neurons. Such models have been used by others as well \cite{GK}. By adjusting this firing rate function, such models can reproduce experimental recordings of neuronal activity (see for instance \cite{berry}). Here, we adjust it for a different purpose, namely to construct a version in which the occurrence of periodic oscillations is tractable.

We consider a population of neurons described by the probability density $n(s,t)$ of finding a neuron in 'state' $s$ at time $t$ 
where $s$ represents the time elapsed since the last discharge. We assume here that the network is homogenous, that is to say that each neuron in the network has the same dynamic. Moreover, we assume that all the neurons are excitatory. According to \cite{cf:pham98a,cf:pham98b} a simple model of the dynamics is given by the age-structured equation 
\begin{equation}
\label{eqneuron}
\left\{ \begin{array}{l}
\frac{\p n(s,t)}{\p t}+ \frac{\p n(s,t)}{\p s}+ p(s,X(t)) \ n(s,t) =0, 
\\[2mm]
N(t):=n(s=0,t)=\int_0^{+\infty} p(s,X(t)) \ n(s,t) ds,
\end{array} \right.
\end{equation}
completed with an initial probability density $n^0(s)$ that satisfies
\beq 
0 \leq n^0(s) \leq 1, \qquad \int_0^\infty n^0(s) ds =1.
\label{as:id}
\eeq  
The coefficient $p(s,X)$ represents the firing rate of neurons in the 'state $s$' and in an environment $X$ resulting from the global neural activity. It is usually small (or even vanishes) during a rest phase that depends on $X$, and increases suddenly afterwards. The density of neurons undergoing a discharge at time $t$ is denoted by $N(t)$ and the boundary condition at $s=0$ means that the neuron re-enters the cycle at 'age $s=0$' after firing. The interactions between the neurons are taken into account through the global neural activity at time $t$, $X(t)$. In the theoretical part of this work we use
\begin{equation}
\label{eq:interaction}
X(t) = N(t)
\end{equation}
which corresponds to instantaneous transmission of neuronal discharge. 
In the numerical part we also take into account synaptic integration, which is modeled by distributed delay function $d$, e.g.,  
$$X(t) =\int_0^t  d(s) N(t-s) ds.$$
Here we also assume that the neurons can discharge earlier when the neural activity is high. We assume that the strength of connectivity is directly related to the sensitivity of the threshold of discharge with respect to the total neural activity.

\vspace{0,5cm}

\noindent These modeling assumptions can be  written on the (nonnegative) coefficients as 
\begin{equation}
\label{as:ps2}
\f{\p}{\p s} p(s,x) \geq 0, \qquad p(s,x)=0 \text{ for } s\in(0,s^*(x)) \quad p(s,x)>0 \text{ for } s>s^*(x), \quad p(s,x){\;}_{\overrightarrow{\; s \rightarrow \infty \; }}\; R(x),
\end{equation}
\begin{equation}
\label{as:ps1}
 \f{\p}{\p x} p(s,x) \geq 0,  \qquad  \quad p(s,x){\;}_{\overrightarrow{\; x \rightarrow \infty \; }}\; 1, \qquad 0<R(x) {\;}_{\overrightarrow{\; x \rightarrow \infty \; }}\; 1 .
\end{equation}
\\
Equations \fer{as:ps2}  and \fer{as:ps1} encompass the key biological properties incorporated into the model. The first is that
following a discharge, a neuron looses its excitability, that is its propensity to discharge in response to a stimulation, 
and recovers it progressively in time. Furthermore, the mean firing rate of the neuron is a function of the stimulus level $x$, 
bounded by $R(x)$. 
The second is that the probability of neuronal firing increases with the intensity of excitatory inputs, 
until some saturation level  is reached. Here again, no further assumptions are made regarding the shape of $p$.

\vspace{0,5cm}

\noindent We will make constant use of two remarkable properties satisfied by solutions to \eqref{eqneuron} (see \cite{Pak_Per_Sal}): a conservation law which expresses the interpretation of $n(t, \cdot)$ as a probability density and an a priori bounds that reflects the normalization $p(s,x) \leq 1$,  
\begin{equation}
\label{conservation}
\int_0^\infty n(s,t) ds = \int_0^\infty n^0(s) ds=1 \qquad \forall t\geq 0,
\end{equation}
\begin{equation}
\label{apb}
0 \leq n(s,t) \leq 1, \qquad  0 \leq N(t) \leq 1, \qquad \forall t\geq 0, \ \forall s >0.
\end{equation}

The theorems and the proofs of the paper are written for a specific choice of the function $p$ which is piecewise constant. More precisely, $p$ is given by
\begin{equation}
\label{assimple}
\left\{ \begin{array}{l}
 p(s,x) =  \ind{s> \sigma (x)}, \quad \text{for some  } \sigma \in Lip(\mathbb{R}^+), \quad  \sigma'(x) \leq 0,
 \\[2mm]
 \sigma(0)=\sigma ^+, \qquad 0< \sigma ^- \leq  \sigma (x) \leq \sigma ^+ , \ \forall x \geq 0.
\end{array}\right. 
\end{equation}
Hence, at each time $t$, the neuron population undergoes two different dynamics as they are or not in their refractory state.
\begin{itemize}
\item  During the refractory period $ (s \leq \sigma(N(t))$, the population density satisfies the free transport equation. At each time $t$, the flux of neurons which enter in their refractory period is given by $N(t)=n(s=0,t)$.
\item After the refractory period,  $(s > \sigma(N(t))$, the solution satisfies the transport equation with an exponential decay due to the neurons discharge. At each time $t$,  neurons enter this second stage where they begin to discharge. As we will see it later, the corresponding flux is given by $n(\sigma(N(t)), t)(1- N'(t) \sigma '(N(t)))$.
\end{itemize}
Therefore, at each time  $t$, the function  $\sigma(N)$ defines the threshold time elapsed when the neurons begin to discharge. The derivative $|\sigma '|$ measures the sensitivity of neurons to the global neural activity. Here our aim is to understand 
the link between  
 $|\sigma'|$ and the appearance of rhythmic and synchronous activity in neural assemblies.  Hence the quantity  that we will essentially study 
in this article is not the function $n(s,t)$, but the density of neurons which discharge at each time 
 $t$, i.e., $N(t)$.  More precisely, we are going to focus on two dynamics.
\begin{itemize}
\item The case where $N$ converges asymptotically to a constant: this dynamic corresponds to a total desynchronization of neuronal discharge
 in the network.
\item The case where $N$ is a periodic function:  this dynamic corresponds to the appearance of rhythmic and synchronous activity, where the higher synchronization of neurons is given for the times   $t$ where $N$ takes its bigger value.
\end{itemize}

\section{Dynamics of ensembles of disconnected neurons  ($ \f{d}{dx} \sigma (x) =0$)}
\label{sec:linear}

To explain our method, we begin with the linear case, when neurons are firing independently from one another. We introduce a tool to prove desynchronization of the population. We first revisit and
 improve the results obtained  in \cite{Pak_Per_Sal} with using a mere  $L^{1}$ norm and thus we remove the exponential weights on the initial data. Secondly, we propose a different analysis, by studying directly the firing rate function $N(t)$. We prove that, as $t$ tends to infinity, $N(t)$ converges to a constant $N^{*}$, oscillating around this constant.

 We assume that $\sigma$ is a small enough constant and thus the equation on $n$ is linear given by
\begin{equation}
\label{eqlin}
\left\{ \begin{array}{l}
\frac{\p n(s,t)}{\p t}+ \frac{\p n(s,t)}{\p s}+ p(s) \ n(s,t) =0, 
\\[2mm]
N(t):=n(s=0,t)=\int_0^{+\infty} p(s) \ n(s,t) ds,
\end{array} \right.
\end{equation}
\begin{equation}
\label{eqlinp}
p(s)=\ind{s \geq \sigma}, \qquad 0 < \sigma<1.
\end{equation}

\subsection{Long time behavior of the solution $n(s,t)$}

\begin{thm}\label{th1}
With assumptions \eqref{as:id} and \eqref{eqlinp},  the solution of the equation (\ref{eqlin}) converges with an exponential rate to the stationary state 
\begin{equation}
\label{eqlinststst}
A(s)=  N^* e^{- \int_{0}^{s} p(s')ds'}, \qquad \text{with } \; N^* = \frac{1}{\sigma +1} .
\end{equation}
 More precisely, there exist $\mu >0$ and $C>0$  such that the following estimate holds
$$  
\int_{0}^{+\infty} |n(s,t)-A(s)|ds \leq C e^{-\mu t} \int_{0}^{+\infty}|n(s,0)-A(s)|ds.
$$
\end{thm}

Several variants of this result are known on bounded intervals (see \cite{MD}  for instance). The interest here is to work on the half-line and in the functional space $L^1$ which is natural in the present context.
\\

\noindent  {\bf Proof of Theorem  \ref{th1}.}
\\
\noindent {\bf Splitting of the solution in two terms.}  The idea of the proof consists in splitting the solution in two parts and set 
$$
n^{0}= n^{0}_{1}+ n^{0}_{2}, \qquad n(s,t)=n_1(s,t)+n_2(s,t).
$$
$\bullet$ The  first initial term $n^{0}_{1}$ takes the value $N^*$ on an interval strictly bigger than  $[0, \sigma]$ and  its support is not necessarily finite.  
\\
$\bullet$ The second  term $n^{0}_{2}$  has zero average and with bounded support, independent of $n^0$ (we can choose $[0,3]$ for instance) and thus we can apply our earlier theory in  \cite{Pak_Per_Sal}.
\\

To build this decomposition, we notice that the stationary solution  $A$ in \eqref{eqlinststst} satisfies
$$
A(s)  \equiv N^{*}  \hbox{ on } [0, \sigma] \qquad  \hbox{ and } \quad  \int_{ \sigma}^{+ \infty} A(s)ds= N^{*}.
$$
Then, we split the initial data as follows
$$
n^{0}_{1}= N^{*}\mathbb{I}_{[0,\sigma+ \frac{1}{2}]} + n^{0}\mathbb{I}_{[\sigma+\frac{1}{2}, +\infty]} - N^{*} \mathbb{I}_{[\sigma+\frac{1}{2},2\sigma +1 ]} 
+ \mathbb{I}_{[\sigma+\frac{1}{2}, \sigma+ \frac{1}{2} + \widetilde{\sigma}(n^{0})]}
$$
with
$$
\widetilde{\sigma}(n^{0})=  \int_{0}^{\sigma+ \frac{1}{2}}n^{0}(s)ds.
$$
We will make use of properties of  $n^{0}_{1}$ and  $n^{0}_{2}$, namely
\\
$\bullet$ $\int_{0}^{+ \infty}n^{0}_{1}(s)ds=1$ and thus  $\int_{0}^{+ \infty} n^{0}_{2}(s)ds=0$.
\\
$\bullet$  The support of $n^{0}_{2}$ is contained in  $[0,3]$ because $\sigma <1$ and $ \widetilde{\sigma}(n^{0}) \leq 1$. 
\\
$\bullet$  $n^{0}_{1}$ is equal to  $N^{*}$ on an interval strictly bigger than $[0,\sigma]$.
\\
$\bullet$ $n(s,t)-A(s)= n_{2}(s,t) + [n_{1}(s,t)-A(s)]$.

\vspace{0,5cm}

\noindent {\bf Time decay of $n_{2}$.}
We recall the following Proposition proved in  \cite{Pak_Per_Sal}, using an entropy method \cite{MMP,Pe}
\begin{prop} \label{entropie}
There exist $\mu >0$ and a continuous function $\psi$, uniformly bounded from below by a positive constant, such that for all bounded initial data $n^{0}$  with $\int_{0}^{+\infty}n^{0}(s)ds=0$, the estimate holds:
$$
\int_{0}^{+ \infty} |n(s,t)| \psi(s) ds \leq e^{-  \mu t} \int_{0}^{+ \infty}  |n^{0}(s)|  \psi (s) ds.
$$
\end{prop}
We can apply Proposition \ref{entropie} to the solution $n_{2}$ and get
$$  
\int_{0}^{+ \infty} |n_{2}(s,t)| \psi(s) ds \leq  e^{-  \mu t} \int_{0}^{+ \infty}  \psi (s) |n^{0}_{2}(s)|ds.
$$
Using the fact that the support on $n^{0}_{2}$ is contained in $[0,3]$, we obtain, with $C=  \frac{\sup_{[0,3]} \psi(s)}{\inf_{\mathbb{R}} \psi(s)}$
\begin{equation}\label{P1}
\int_{0}^{+ \infty} |n_{2}(s,t)| ds \leq C  e^{-  \mu t}\int_{0}^{+ \infty}  |n^{0}_{2}(s)|ds.
\end{equation}

\vspace{0,5cm}

\noindent {\bf Estimate of $n_{1}-A$.}  We begin with a Lemma that will be used on $n_1-A$,
\begin{lemme}\label{sta}
Let  $n^{0}$ be such that $\int_{0}^{+\infty}n^{0}(s)ds=0$ and  $n^{0} \equiv 0 $  on $[0,\sigma]$. Then, we have
$$
N(t) \equiv n(s,t) \equiv 0 \quad  \hbox{ for } s\in [0,\sigma], \qquad t\geq 0.
$$
\end{lemme}
{ \bf Proof of Lemma \ref{sta}.}
Using the characteristics, we deduce that for all time $t$ and $s \in [0,\sigma]$
$$
n(s,t)= N(t-s) \qquad \hbox{ where we extend $N$ for}\; t \in [- \sigma,0]  \hbox{ by } \; N(t)=0.
$$
With this extension, the mass conservation gives
\begin{equation}\label{mas}
N(t)+ \int_0^{\sigma} N(t-s)ds=0, \qquad \forall t \geq 0.
\end{equation}
Because $N(t)=0$ for $t\in [-\sigma,0]$, we conclude that for all $\tau >0$
$$
\sup_{0\leq t \leq \tau} |N(t)|�\leq\sup_{0\leq t \leq \tau}  \int_0^{\sigma} |N(t-s)|ds \leq \sigma \sup_{0\leq t \leq \tau}  |N(t)|
$$
which implies (because $\sigma <1$) that $N\equiv 0$, which ends the proof of Lemma \ref{sta}. 
 \hfill $\square$

\medskip

\noindent 
To obtain the exponential decay in time of  $ \|n_{1}-A (t) \|_{L^{1}}$, we are reduced to proving the following Lemma (which we will apply to  $n_{1}-A$ also)

\begin{lemme} \label{sseuil}
Let $n^{0}  \in L^{1}$ with $\int_{0}^{+\infty}n^{0}(s)ds=0$ be a function such that $n^{0}(s)=0$ for $s\in [0, \sigma]$ 
 and let $n$ be the solution of the equation
\begin{equation}
\label{eqss}
\left\{ \begin{array}{l}
\frac{\p n(s,t)}{\p t}+ \frac{\p n(s,t)}{\p s}+  \ n(s,t) =0, 
\\[2mm]
N(t):=n(s=0,t)=\int_\sigma^{+\infty} \ n(s,t) ds.
\end{array} \right.
\end{equation}
Then the following estimate holds
\begin{equation}
\label{eqssexpd}
\int_{0}^{+ \infty} |n(s,t)|ds  \leq e^{-t}  \int_{0}^{+ \infty}|n^{0}(s)| ds.
\end{equation}
\end{lemme}
 
\vspace{0,5cm}

\noindent {\bf Proof of Lemma \ref{sseuil}.}   
Since the  function $|n|$  is solution of 
$$
\partial_{t}|n|+ \partial_{s}|n|+|n|=0, 
$$
using Lemma \ref{sta}, we know that  $n(s,t) \equiv 0$ for $s\in [0,\sigma]$. Therefore 
we may integrate with respect to the variable $s\in [0,+\infty[$, we deduce that 
$$
\frac{d}{dt}\int_{0}^{+\infty}|n|(s,t)ds=- \int_{0}^{+\infty}|n|(s,t)ds=- \int_{\sigma}^{+\infty}|n|(s,t)ds
$$
from which we conclude the inequality \eqref{eqssexpd} and Lemma  \ref{sseuil} is proved.  
\hfill $ \square$

\medskip

With Lemma \ref{sseuil} we conclude that 
$$ 
\int_{0}^{+ \infty} |n_{1}(s,t)-A(s)| ds \leq e^{-t}  \int_{0}^{\infty} |n^0_{1}(s)-A(s)| ds .
$$
Combined with the  estimate (\ref{P1}), the proof of Theorem \ref{th1} is complete.  
\hfill $ \square$

\subsection{Oscillatory relaxation of the global activity $N$}
\label{sec:N}

As a preparation to the nonlinear case,  we now explain another method to prove asymptotical total desynchronization. This method relies on a direct study of the function $N(\cdot)$. 

We mention that, as we will see it later, we are in the simple situation where at each time  $t \geq \sigma$, the flux of neurons which enter in the area  of discharge is given by 
$$ 
n(\sigma,t) [1-\frac{dN(t)}{dt}\sigma '\big(N(t)\big)] = N(t-\sigma).
$$

More precisely the following Proposition holds:
\begin{prop}\label{prop1}
With the assumptions of Theorem \ref{th1}, the solution $N$ converges exponentially with rate $\ln \sigma$ to the constant $N^{*}= \frac{1}{1+\sigma}\cdotp$ Moreover, 
$N(t)$ oscillates around of $N^{*}$; for any interval $I$ of size  $\sigma$
$$
\exists t_0 \in I  \; \hbox{ such that  } \quad N(t_0) = N^{*} .
$$
\end{prop}

\noindent {\bf Proof of Proposition \ref{prop1}.}
Using the characteristics for equation \eqref{eqlin} and mass conservation, we obtain 
\beq
N(t)+ \int_{t-\sigma}^{t}N(s)ds=1, \quad  \forall t \geq \sigma,
\label{relation_lin}
\eeq
which we write as 
$$
N(t)-N^{*}= \int_{t-\sigma}^{t} [N^{*}  - N(s)] ds.
$$
Using that $0 \leq N^* \leq 1$, $0 \leq  N \leq 1$ in the right hand side, we obtain that 
$$
|N(t)-N^{*}| \leq \sigma, \quad  \forall t \geq \sigma.
$$
We may iterate the above argument and deduce that the following estimate holds
$$
|N(t)-N^{*}| \leq \sigma^{n}, \quad  \forall t \geq n \sigma,
$$
which ends the proof of the exponential convergence of $N$ to $N^{*}$ because $ \sigma <1$.
\\

 It is also easy to prove that the solution oscillates around $N^{*}$. Let $I= (t_{0}, t_{0}+ \sigma)$ be an interval 
of size $\sigma$. The relation \eqref{relation_lin} implies that
$$
N(t_{0}+\sigma)+ \int_{I} N(s)ds=1,
$$
which contradicts the equality $N^{*}+ \int_{t-\sigma}^{t}N^{*} ds=1$ 
if $N(t)> N^{*}$ or $N(t) < N^{*}$ on $I$. This completes the proof of  Proposition  \ref{prop1}. 
\hfill $ \square$

\section{ Asynchronous activity for connected neuronal assemblies  ($\frac{d \sigma}{dx} \neq 0$)}
\label{sec:desync}

We come back to the nonlinear situation when the threshold depends on the global neural activity to study relaxation properties. 
In \cite{Pak_Per_Sal}  it was proved that both for weakly and highly connected networks, the model undergoes relaxation to the unique steady state. But the precise conditions in there are difficult to follow and are certainly not optimal. Our purpose here is to come to precise (and somehow sharp) conditions leading to global relaxation. We introduce a different strategy, following the method of section \ref{sec:N}. 

Our first condition is that 
for some  $0 \leq m <1$ we have 
\begin{equation}\label{Z1}
0 \leq - \sigma ' (x) \leq m <1 .
\end{equation}
It implies that for any probability density $0\leq n(\cdot) \leq 1$, there is a unique $N$ such that
$$
N=\int_{\sg(N)}^{\infty} n(s)ds, \qquad 0<N<1,
$$
because $N\mapsto \int_{\sg(N)}^{\infty} n(s)ds$ is a contraction and we can apply the Banach-Picard Theorem. In particular the boundary condition \eqref{eqneuron} defines a unique $N(t)$, while the periodic solutions we build later correspond to two possible choices of $N$, hence discontinuities.

\subsection{Some observations}
 
The conservation law \eqref{conservation} can be written
\begin{equation}
\label{eq:basicN}
1= N(t)+ \int_{0}^{\sigma(N(t))}n(s,t)ds.
\end{equation}
The following lemma ensures that $N$ is smooth enough for our purpose
\begin{lemme}\label{lipsch}
 Assume that condition \eqref{Z1}  holds. Then, the function $N$ is Lipschitz continuous.
\end{lemme}

\noindent {\bf Proof of Lemma \ref{lipsch}.} 
Consider two times $t_1$, $t_2$. The relation \eqref{eq:basicN} gives for $i=1,2$
$$
1=N(t_{i})+ \int_{0}^{\sigma(N(t_{i}))}n(s,t_{i})ds,
$$
hence
$$
N(t_{1})-N(t_{2})= \int_{0}^{\sigma(N(t_{1}))}[n(s,t_{2})-n(s,t_{1})]ds+ \int_{\sigma(N(t_{1}))}^{\sigma(N(t_{2}))}n(s,t_{2})ds.
$$
Using the fact that  $0 \leq n\leq1$ and condition \eqref{Z1}, we find that
\begin{equation}\label{es1}
\left|\int_{\sigma(N(t_{1}))}^{\sigma(N(t_{2}))}n(s,t_{2})ds \right| \leq m |N(t_{1})-N(t_{2})|.
\end{equation}
On the other hand, using the PDE on $n$, we have
$$\begin{array}{rl}
\int_{0}^{\sigma(N(t_{1}))}  [n(s,t_{2}) - &n(s,t_{1}) ] ds= \int_{t_{1}}^{t_{2}}\int_{0}^{\sigma(N(t_{1}))} \partial_{t} n(s,t)ds dt
\\ \\
& = \int_{t_{2}}^{t_{2}}\int_{0}^{\sigma(N(t_{1}))} \partial_{s} n(s,t)+ p(s,N(t))n(s,t)ds dt
\\ \\
&= \int_{t_{2}}^{t_{1}}\int_{0}^{\sigma(N(t_{1}))}     p(s,N(t))n(s,t)ds dt + \int_{t_{1}}^{t_{2}} [n(  \sigma(N(t)),t)-N(t) ] dt.
\end{array}
$$
We obtain that there exists a constant $C$ such that
\begin{equation} \label{es2}
\left|\int_{0}^{\sigma(N(t_{1}))}[ n(s,t_{2})-n(s,t_{1})] ds \right| \leq C |t_{1}-t_{2}|.
\end{equation}
Combining estimates (\ref{es1}) and (\ref{es2}) and the fact that $m<1$, we deduce that there exists a constant $C(m)$ such that 
$$
|N(t_{1})-N(t_{2})| \leq C(m)|t_{1}-t_{2}|
$$
which ends the proof of Lemma \ref{lipsch}. 
\hfill $\square$
  
 \vspace{0,5cm}
  
We may now differentiate a.e. the relation  (\ref{eq:basicN})  and find 
\begin{equation}
\label{eq:fluxn}
\f{dN(t)}{dt}  \big[1+ \sigma '\big(N(t)\big) n\big(\sigma(N(t)),t\big) \big]= -N(t)+ n\big(\sigma(N(t)),t\big).
\end{equation}
This equation is much more complicated to use than in the case when $\sigma$ is constant. It can be interpreted as follows. At each time  $t$,
\\
$\bullet$ the boundary condition tells us that the flux of neurons which enter in their refractory state is given by $N(t)$.
\\
$\bullet$ the flux  of neurons which enter in the state of possible discharge is given by 
$$
n\big( \sigma(N(t)),t \big)\;  [1-\frac{dN(t)}{dt}\sigma '\big(N(t)\big)] .
$$
This quantity is not as easy to control as in the linear case. For instance  
it can take nonpositive values if the velocity of the threshold given by $\f{dN(t)}{dt} \sigma'(N)$ gets bigger than the chronological speed.

\subsection{A condition for asymptotic total desynchronization}

In this section, using assumption \eqref{Z1}, we give a sufficient condition on the function $\sigma$ and the maximal refractory duration 
$$
\sigma^{+}= \sup_{\RR^+}\sigma(x)= \sg(0)
$$
that implies asymptotical total desynchronization of the neural activity. 
\begin{thm}\label{TH1}
Assume \eqref{Z1} and that, with $\bar N <1$ the unique steady state defined in \eqref{barn} below,   
\begin{equation}\label{Z2}
\sg(x) \leq \sigma ^{+} < 1-  m \bar N,
\end{equation}
Then, we have $\f{dN(t)}{dt}  \sigma'(N(t))<1$ and with an exponential rate,
$$
 \lim_{t \to +\infty} N(t) = \bar{N}.
$$
\end{thm}

In section \ref{sec:periodic}, we give explicit examples where the equality $\f{dN(t)}{dt}  \sigma'(N(t))=1$ holds and $N(t)$ can be periodic. 
\\

\noindent {\bf Proof of Theorem \ref{TH1}.}
Assumption \eqref{Z1} and the a priori estimate $0\leq n(s,t) \leq 1$ imply that 
$$
1 \geq g(t) := 1+ \sigma '\big(N(t)\big) n\big(\sigma(N(t)),t\big) \geq 1-m >0.
$$
Hence, using the relation \eqref{eq:fluxn}, we conclude  
$$
\f{dN(t)}{dt}  =  \frac{-N(t)+ n(\sigma(N(t)),t)}{g(t)}.
$$
 
We decompose the end of the proof as follows
\\
$\bullet$  In a first step, we prove that condition (\ref{Z1}) implies that the threshold velocity 
is always strictly less than  $1$ as soon $t \geq \sigma ^{+}$,  i.e., the following estimate holds
\begin{equation}\label{vitseuil}
\f{dN(t)}{dt} \sigma'\big(N(t)\big) \leq  m, \qquad  \forall  t \geq \sigma^{+},
\end{equation}
a property (\ref{vitseuil}) that allows us to simplify the equation on $N$.
\\
$\bullet$ In a second step,  this simplified equation on  $N$ can be used to show that  $N$ converges, when $t \to + \infty$. This is based on an  iterative process similar to that used in section \ref{sec:linear}. 
\\

As a first step we prove the

\begin{lemme}\label{vs}
With condition  (\ref{Z1}), the estimate (\ref{vitseuil}) holds.
\end{lemme}

\noindent {\bf Proof of Lemma \ref{vs}}.
Using the  relation \eqref{eq:fluxn} on $N$, we deduce that, with the notation $\wt n (t)= n(\sigma(N(t)),t)) \geq 0$, we have (still because $\sigma' \leq 0$)
$$\begin{array}{rl}
\f{dN(t)}{dt}  \sigma'\big(N(t)\big) &= \f{ \sigma'\big(N(t)\big) [\wt n (t) - N(t)]}{1+ \sigma'\big(N(t)\big)\;  \wt n (t)} 
\\[3mm]
&= 1+   \f{-1- \sigma'\big(N(t)\big) N(t)}{1+ \sigma'\big(N(t)\big)\;  \wt n (t)}
\\ [3mm]
& \leq 1 + (-1+ m)=  m.
\end{array}$$
\hfill $\square$

\medskip

Next, we simplify the equation on   $N$ using that the threshold velocity is strictly less than $1$.

\begin{lemme}\label{masse}
Assume that estimate  (\ref{vitseuil}) holds.  Then, we have
\begin{equation} \label{eq:integralform}
N(t)+ \int_{t-\sigma (N(t))}^{t}N(s)ds=1 \qquad \forall t\geq \sigma ^{+}.
\end{equation}
\end{lemme}

\noindent{\bf Proof of Lemma \ref{masse}.}
Using the characteristics, we obtain that
$$
n(s,t)= N(t-s)e^{- \int_{0}^{s}p(s',N(t+s'-s))ds'}  \quad  \hbox{for }\ s \leq \sigma(N(t)). 
$$
Take $s  \in [0, \sigma(N(t))]$ and let us prove that
$$
n(s,t)= N(t-s).
$$
To do this, it enough to show that in the above integral $p(s',N(t+s'-s))=0$, that is 
\begin{equation}\label{mes}
  \forall s' \in [0,s) \ \hbox{ we have  } \  s' < \sigma(N(s'+t-s)).
\end{equation}
We set
$$
f(s'):=s'- \sigma(N(s'+t-s)).
$$
We have $f(s)\leq 0$  because $s \leq \sigma(N(t))$. Moreover, from  (\ref{vitseuil}), we deduce that
$$
f'(s')= 1- \f{dN(s'+t-s)}{dt}  \sigma'(N(s'+t-s)>0,
$$
which proves (\ref{mes}). We conclude the proof of Lemma \ref{masse} using \eqref{eq:basicN} which now reads
$$
1= N(t)+ \int_{0}^{\sigma(N(t))}n(s,t)ds=N(t)+ \int_{0}^{\sigma(N(t))}N(t-s)ds.
$$
\hfill $\square$

\bigskip

 We can now conclude the proof of Theorem  \ref{TH1}. We notice there is a unique steady state solution $0 < {\overline N} <1$, that is a solution  to 
\begin{equation}\label{barn}
 {\overline N}  \; \big(1+\sigma ( {\overline N} ) \big) =1.
\end{equation}
This is because the mapping $N \mapsto 1- N \sigma ( N )$ is a contraction as soon as $0 < \sigma(\cdot) <1$ and  $-1 < \sigma'(\cdot) \leq 0$ and thus it has a unique fixed point by the Banach-Picard theorem. 

Then, for all $t \geq \sigma ^{+}$ we rewrite \eqref{eq:integralform} as 
$$
N(t)- {\overline N}+ \int_{t-\sigma (N(t))}^{t}[N(s)-{\overline N}]ds+{\overline N} [\sigma (N(t))- \sigma ({\overline N})]=0.
$$
This proves that for $t \geq \sigma ^{+}$
$$
| N(t)- {\overline N} | -{\overline N} | \sigma (N(t))- \sigma ({\overline N})| \leq \int_{t-\sigma (N(t))}^{t} | N(s)-{\overline N} | ds,
$$
$$
 | N(t)- {\overline N}  \leq \f{\sigma^+}{1-m{\overline N} } \; \sup_{t-\sigma (N(t))\leq s \leq t} \; | N(s)-{\overline N} |.
$$
Since $0 \leq {\overline N}, N(\cdot) \leq 1$, we find $ | N(t)- {\overline N} | \leq \f{\sigma^+}{1-m{\overline N} } $.

We may iterate and for $t \geq n \sigma ^{+}$
$$
 | N(t)- {\overline N} | \leq  \left[ \f{\sigma^+}{1-m{\overline N}}\right]^n.
$$
Using condition \eqref{Z2}, this ends the proof of Theorem \ref{TH1}. 
\hfill $\square$

\section{Synchronous periodic oscillations in neuronal assemblies}
\label{sec:periodic}

In our previous study \cite{Pak_Per_Sal}, we have presented numerical examples of periodic solutions corresponding to synchronous activity that were obtained in numerical simulations. Here, for a particular class of functions $\sigma(\cdot)$ we construct analytically such periodic solutions. To do this, we impose a priori that, on a given open interval $I$,  or a union of two intervals, neurons do not leave their refractory state, namely
\begin{equation}\label{perel0}
\f{d}{dt}  \sigma \big(N(t) \big)=1 \qquad \forall t \in I .
\end{equation}
This can only last until saturation occurs, leading suddenly to a massive discharge of neurons characterized by a discontinuity of $N(t)$. 
\\

For $\alpha >0$, we define the two functions 
\begin{equation} \label{def_N}
0< N^{-}(\alpha) := \frac{1}{2e^{\alpha}-1}\; < \, N^{+}(\alpha):= \frac{e^{\alpha}}{2e^{\alpha}-1}\; < \; 1 ,
\end{equation}
and we choose the Lipschitz continuous discharge threshold $ \sigma$ as 
\begin{equation} \label{def_sigma}
\sigma(x) =   
\left\{ \begin{array}{ll}
2  \alpha   & \hbox{ on } \; [0, N^{-}(\alpha)] ,
\\
2  \alpha -  \ln(x) +  \ln(  N^{-}(\alpha)) & \hbox{ on } \; [N^{-}(\alpha),N^{+}(\alpha)] ,
\\
\alpha  & \hbox{ on } \;  [N^{+}(\alpha), \infty) .
\end{array} \right.
\end{equation}

Motivated by numerical results that we show in section \ref{sec:num}, we are going to build three distinct classes of periodic solutions. The first two are simpler than the third, but they seem unstable because only the third class is accessible by our numerical simulations. As we will see, there are several common features in our constructions. One of them is to postulate a form of $N(t)$ and solve the linear transport equation 
\begin{equation}\label{eqperio}
 \left\{ \begin{array}{l}
\frac{\p n(s,t)}{\p t}+ \frac{\p n(s,t)}{\p s}+ p(s,N(t)) \ n(s,t) =0,  \qquad t \in \mathbb{R}, \; s\geq 0, 
\\[2mm]
n(s=0,t) = N(t).
\end{array} \right.
\end{equation}
Because it is periodic in $t$, there is no need of initial data and the method of characteristics gives a solution. The three classes depend on suitable choices of $N(t)$. 

\subsection{Periodic solutions with one discontinuity}
\label{sec:periodic1}

The first class is characterized by a periodic activity $N(t)$ with a single discontinuity per period. 
   
\begin{thm} \label{perio}
For $\alpha >0$, define the discontinuous periodic function $N(t)$ of period $ \alpha$ as 
$$
 N(t)= N^{+}(\alpha) e^{-t} \;  \hbox{ for } \; t  \in I:= (0,  \alpha) .  
$$
Then the solution to the linear renewal equation \eqref{eqperio} is also a solution of the equation (\ref{eqneuron}). In other words, we have for all $t\in \mathbb{R}$
 $$
 N(t) = \int_{\sg(N(t))}^{\infty} n(s,t)ds \qquad  \text{and} \qquad \int_0^\infty n(s,t) ds =1.
 $$
\end{thm}
 
Some comments about this Theorem are in order
\\
$\bullet$ Since $N(\alpha^+)=N^-(\alpha)$,  $N(t)\in ( N^{-}(\alpha), N^{+}(\alpha))$ is a  sawtooth function; it decreases exponentially from $N^{+}(\alpha)$ to $N^{-}(\alpha)$. When the discharge rate of neurons is at its minimal value $N^{-}(\alpha)$, a peak of discharge occurs in the network which is characterized by a jump of the function $N$ from $N^{-}(\alpha)$ to $N^{+}(\alpha)$. Consequently the precise shape of $\sigma$ outside the interval $[N^{-}(\alpha), N^{+}(\alpha)]$ does not play a role here.
\\
$\bullet$ We can compute that, for these solutions, the assumption \eqref{Z1} does not hold. Indeed
$$
\sup |\sigma'(N(t))|= 1/N^-(\alpha)= 2e^{\alpha}-1 >1.
$$
\\
$\bullet$ The bigger $\alpha$ is, the bigger is the difference $N^{+}(\alpha)-N^{-}(\alpha)$. In particular, $N^{+}(0)-N^{-}(0)=0$ and
$$
\lim_{\alpha \to +\infty} N^{+}(\alpha)=\frac{1}{2} \quad \hbox{ and } \quad  \lim_{\alpha \to +\infty} N^{-}(\alpha)=0.
$$
The synchronization phenomenon is more evident for large values of $\alpha$. Hence we can hope that sustained oscillations still occur for  more general models when $\alpha$ is large
enough. For instance, numerical simulations indicate that this choice of the function $\sigma$ also yields periodic solutions when the global neural activity is defined with a delay as in \cite{Pak_Per_Sal}.
\\

\noindent{\bf Proof of Theorem  \ref{perio}.} 
%
To begin with, we notice that \eqref{perel0} holds true with our choice. Then we argue in two steps. First we establish a form of the boundary condition, then the integral condition.
\\

\noindent{\bf First step.} We notice that, for $t\in I$, the discontinuities of $n(s,t)$ lies on $\{ s=k \alpha + t\}$. Therefore, taking into account \eqref{perel0},  we calculate, for $t\in I$,
$$
\begin{array}{rl}
\f{d}{dt} \int_0^{\sigma(N(t))} n(s,t) ds & = \f{d}{dt}  \sigma (N)(t) \; n\big(\sigma(N(t))^-,t\big)+ \int_0^{\sigma(N(t))} \f{\partial}{\partial t} n(s,t) ds 
\\[3mm]
& =  n\big(\sigma(N(t))^-,t\big) - \int_0^{\sigma(N(t))} \f{\partial}{\partial s} n(s,t) ds 
\\[3mm]
&= n(0,t)=N(t)=-N'(t).
\end{array}
$$
We may now write, by periodicity, that for all $t\in \mathbb{R}$,
\begin{equation}\label{perrel1}
 \int_0^{\sigma(N(t))}  n(s,t) ds + N(t) := Q_1=1.
\end{equation}
Indeed, in order to evaluate the constant $Q_1$, we can use the method of characteristics and compute this quantity at $t=0^+$
$$
\int_0^{\sigma(N(0^+))} n(s,0) ds = \int_0^{\sigma(N^+(\alpha))} N(-s) ds =N^+(\alpha) \int_0^{\alpha} e^{s-\alpha} ds=N^+(\alpha) (1-e^{-\alpha}) .
$$
After evaluating this expression, we arrive at 
$$
\int_0^{\sigma(N(0^+))} n(s,0) ds =  1-N^+(\alpha)  , \qquad Q_1=1.
$$

\noindent{\bf Second  step.} Finally, integrating equation \eqref{eqperio}, we compute for $t\in \mathbb{R}$ (the second equality uses  \eqref{perrel1})
$$
\f{d}{dt} \int_0^{\infty} n(s,t) ds = N(t)-  \int_{\sigma(N(t))}^\infty  n(s,t) ds =1- \int_0^{\infty} n(s,t) ds.
$$ 
For a periodic function this proves that $ \int_0^{\infty} n(s,t) ds=1$ and the proof of Theorem \ref{perio} is complete.
\hfill $\square$

\subsection{Periodic solutions with two discontinuities}
\label{sec:periodic2}

It is possible to build a family of more elaborate periodic solutions by allowing for two discontinuities of $N(t)$ during a period. The construction uses an additional parameter 
\beq
p \in (0,\alpha), \qquad N_p(\al) = e^{p} \; N_-(\alpha) \in \big(N_{-}(\alpha),N_{+}(\alpha) \big),
\label{eq:functiongNp}
\eeq
and the single root $\gamma(p) \in ]p+\alpha,2 \alpha]$ to the equation $g(\gamma)=0 $ where 
\beq
g(\gamma)= e^{-\alpha}[e^{\al-p}-1][e^{\gamma-\alpha}-\gamma+1]
+e^{-\alpha}[-2\alpha+( \alpha +p)e^{\al-p}]+e^{p-\al}- 1.
\label{eq:functiong}
\eeq
\begin{thm}\label{perio2}
For $p \in (0,\alpha)$, define the $\gamma(p) $ periodic function by
$$N(t)= 
\left\{\begin{array}{ll}
N_{+}(\alpha) \; e^{-t} &   \hbox{ for } \; t  \in  [0,\alpha] ,
\\
N_p(\al)  \; e^{-t+\alpha} &  \hbox{ for } \; t  \in  (\alpha,p+\alpha], \\
N_{+}(\alpha)\; e^{-t+\alpha} [t- \gamma(p)+e^{\gamma(p)-\alpha}] &  \hbox{ for } \; t  \in  (p+\alpha,\gamma(p)).
\end{array}\right.
$$
Then the solution to the linear renewal equation \eqref{eqperio} is also a solution of the equation (\ref{eqneuron}). In other words, we have for all $t\in \mathbb{R}$
 $$
 N(t) = \int_{\sg(N(t))}^{\infty} n(s,t)ds \qquad  \text{and} \qquad \int_0^\infty n(s,t) ds =1.
 $$
\end{thm}

The single peak solution in section \ref{sec:periodic1} corresponds to the limiting case $p=0$, $\gamma(p)=\alpha$.
\\

\noindent {\bf Proof.}
As in the proof of Theorem \ref{perio}, we have for $t \in (0, \alpha)$ and $t \in (\alpha,p+\alpha)$
\begin{eqnarray}
\frac{d}{dt} \int_{0}^{\sigma(N(t))}n(s,t)ds & = &N' \sigma'(N(t))n(\sigma(N(t)),t) + \int_{0}^{\sigma(N(t))} \partial_{t} n(s,t)ds
\nonumber \\ \nonumber 
&=& n(\sigma(N(t)),t)-\int_{0}^{\sigma(N(t))} \partial_{s} n(s,t)ds 
\nonumber \\ \nonumber 
&=& n(0,t)=N(t)=-N'(t).
\end{eqnarray}
For $t \in ( \alpha+p, \gamma(p))$, we have $N(t) \geq N_{+}(\alpha)$. Therefore $\sigma(N(t))= \alpha$ and we calculate
\begin{eqnarray}
\frac{d}{dt} \int_{0}^{\sigma(N(t))}n(s,t)ds & = & N' \sigma'(N(t))n(\sigma(N(t)),t) + \int_{0}^{\sigma(N(t))} \partial_{t} n(s,t)ds
\nonumber \\ \nonumber 
&=& -\int_{0}^{\sigma(N(t))} \partial_{s} n(s,t)ds \nonumber \\ \nonumber 
&=& n(0,t)-n(\alpha,t)=N(t)-N(t-\alpha)
\nonumber \\ \nonumber 
&=& N_{+}(\alpha) e^{-t+\alpha} [t-\gamma(p)+e^{\gamma(p)}]- N_{+}(\alpha) e^{-t+\alpha} =-N'(t).
\end{eqnarray}
By periodicity, we conclude that there is a constant $m(p)$ such that for all $t\in \RR$
$$
\int_{0}^{\sigma(N(t))}n(s,t)ds+N(t)=m(p).
$$
Furthermore, following the second step in the proof of Theorem  \ref{perio}, we have for all $t \geq 0$
$$
\int_{0}^{+\infty} n(s,t)ds= m(p).
$$
It remains to check that $m(p)=1$ and it is enough to prove that
$$
\int_{0}^{\sigma(N(0^{+}))}n(s,0^{+})ds+N(0) =1.
$$
That is to say 
$$
\frac{1}{N_{+}(\alpha)}\int_{\alpha+p}^{\gamma(p)}N(s)ds+e^{p-\al}= 1.
$$
We have
$$ 
\frac{1}{N_{+}(\alpha)}\int_{\alpha+p}^{\gamma(p)}N(s)ds= e^{-\alpha}[ e^{\al-p} -1][e^{\gamma(p)-\alpha}-\gamma(p)+1]
+e^{-\alpha}[-2\alpha+(\alpha +p)e^{\al-p}].
$$
With the definition in \eqref{eq:functiong}, we deduce that we must have $g(\gamma(p))=0$. 
\\

It remains to prove that there exists a unique value $\gamma(p) \in (\alpha+p,2\alpha)$ such that $g(\gamma(p))=0$. 
To prove this, we observe that 
$$
g(\alpha+p)= e^{-\alpha}[e^{\al-p} + p- \al -1] \leq 0,
$$
and that
$$
g(2\alpha)=e^{-\alpha} e^{\al-p}[e^{\alpha}+p-\al+1]-2+ e^{p-\al} -e^{-\alpha} \geq 0.
$$
Since $g(\cdot)$ is increasing, this completes the proof of Theorem \ref{perio2}.
\hfill $\square$

\subsection{Periodic solutions with two discontinuities and a flat state}
\label{sec:periodic3}

We are now ready to present the third class of periodic solutions which seems to be those observed numerically. We use again the notations \eqref{eq:functiongNp} and we need the function
\beq
f(\delta,y):=e^{\alpha}\Big([-2\alpha+1+e^{\alpha}+\delta]e^{-\delta}-e^{-2\alpha}-e^{-\alpha} \Big)+
e^{-\alpha}[\delta-\ln(y)-2\alpha]+y-1.
\label{eq:functionf}
\eeq

\begin{thm}\label{perio1}
For $p\in (0,\al)$ and $\delta(p) \in [p+\alpha, 2\al]$, let $Y(p)$ be such that $f(\delta(p),Y(p))=0$ (we will see that $\delta(p)$ exists), let  $N$ be the $2 \alpha$ periodic function given by
$$N(t)= 
\left\{\begin{array}{ll}
N_{+}(\alpha) \; e^{-t} &   \hbox{ for } \; t  \in  [0,\alpha] ,
\\
N_p(\al)  \; e^{-t+\alpha} &  \hbox{ for } \; t  \in  (\alpha,p+\alpha], 
\\
N_{-}  (\al)    &  \hbox{ for } \; t  \in   (p+\alpha,\delta(p)],
\\
N_{+}(\alpha)\; e^{-t+\alpha}[t-2 \alpha+e^{\alpha}] &  \hbox{ for } \; t  \in   (\delta(p),2 \alpha).
\end{array}\right.
$$
Then the solution to the linear renewal equation \eqref{eqperio} is also a solution of the equation (\ref{eqneuron}). In other words, we have  for all $t\in \mathbb{R}$
$$
 N(t) = \int_{\sg(N(t))}^{\infty} n(s,t)ds \qquad  \text{and} \qquad \int_0^\infty n(s,t) ds =1.
$$
\end{thm}

The double peak solution in section \ref{sec:periodic2} corresponds to $\delta= p+\alpha$ and $Y(p)=e^{p-\al}$ (no flat state).
\\

\proof 
As before, we can write  for $t \in (0, \alpha)$ and $t \in (\alpha,p+\alpha)$
\begin{eqnarray}
\frac{d}{dt} \int_{0}^{\sigma(N(t))}n(s,t)ds & = &N' \sigma'(N(t))n(\sigma(N(t)),t)-\int_{0}^{\sigma(N(t))} \partial_{s} n(s,t)ds
\nonumber \\ \nonumber 
&=& n(0,t)=N(t)=-N'(t).
\end{eqnarray}
For $t \in (\alpha+p,\delta(p))$,  $N(t) \equiv N_{-}(\alpha)$ and so $\sigma(N(t))\equiv 2 \alpha$  on $(\alpha+p,\delta(p))$. Hence
for $t \in (\alpha+p,\delta(p))$,  
\begin{eqnarray}
\frac{d}{dt} \int_{0}^{\sigma(N(t))}n(s,t)ds & = & N' \sigma'(N(t))n(\sigma(N(t)),t)-\int_{0}^{\sigma(N(t))} \partial_{s} n(s,t)ds
\nonumber \\ \nonumber 
&=& -\int_{0}^{\sigma(N(t))} \partial_{s} n(s,t)ds \nonumber \\ \nonumber &=& n(0,t)-n(2\alpha,t)=
N(t)-N(t-2\alpha)=0=-N'(t).
\end{eqnarray}
For $t \in ( \delta(p), 2\alpha)$, we have $N(t) \geq N_{+}(\alpha)$. We deduce that $\sigma(N(t))= \alpha$ and that
\begin{eqnarray}
\frac{d}{dt} \int_{0}^{\sigma(N(t))}n(s,t)ds & = & N' \sigma'(N(t))n(\sigma(N(t)),t)-\int_{0}^{\sigma(N(t))} \partial_{s} n(s,t)ds
\nonumber \\ \nonumber 
&=& -\int_{0}^{\sigma(N(t))} \partial_{s} n(s,t)ds \nonumber \\ \nonumber &=& n(0,t)-n(\alpha,t)=N(t)-N(t-\alpha)
\nonumber \\ \nonumber 
&=& N_{+}(\alpha) e^{-t+\alpha} [t-2 \alpha+e^{\alpha}]- N_{+}(\alpha) e^{-t+\alpha} =-N'(t).
\end{eqnarray}

Again, we conclude by periodicity that there is a constant $m(p)$ such that for all $t\in \RR$
\begin{equation}\label{consm}
\int_{0}^{\sigma(N(t))}n(s,t)ds+N(t)=m(p).
\end{equation}
To prove that   for all $t \geq 0$
$$
\int_{0}^{+\infty} n(s,t)ds= m(P),
$$
we use the equality (\ref{consm}),  and we have
$$
\frac{d}{dt}\int_{0}^{+\infty} n(s,t)ds= N(t)- \int_{\sigma(N(t))}^{+\infty} n(s,t)ds=m(p)- \int_{0}^{+\infty}n(s,t)ds.
$$
By periodicity, we deduce that  for all $t \geq 0$ 
$$
\int_{0}^{+\infty} n(s,t)ds=m(p).
$$
Let us assume that the following lemma holds and let us  prove that $m(p)=1$.
\begin{lemme}\label{estf} Recalling the function $f$ in \eqref{eq:functionf}, for all $Y(p) \in [e^{-\alpha},1]$, there exists $\delta(p) \in [p+\alpha,2\alpha]$ such that, 
$$
f(\delta(p),Y(p))=0.
$$
\end{lemme}
To prove that $m(P)=1$, it is enough to check that
$$
\int_{0}^{\sigma(N(0^{+}))}n(s,0^{+})ds+N(0) =1,
$$
that is 
\begin{equation}\label{estmasse}
\int_{\alpha}^{2\alpha}N(s)ds= 1-N_{+}(\alpha).
\end{equation}
We have 
$$
\int_{\alpha}^{p+\alpha} N(s)ds=N_{p}(\alpha) P-N_{-}(\alpha) 
$$
hence, by setting $Y(p)=e^{p-\alpha}$, we rewrite estimate (\ref{estmasse}) as
$$
\frac{1}{N_{+}(\alpha)}\int_{\delta(p)}^{2\alpha}N(s)ds +e^{-\alpha}[\delta(p)-\ln(Y(p))-2\alpha]+Y(p)= 1.
$$
We have
$$
\int_{\delta(p)}^{2 \alpha}N(s)ds  = N_{+}(\alpha)e^{\alpha}\Big([-2\alpha+1+e^{\alpha}+\delta(p)]e^{-\delta(p)} -e^{-2\alpha}-e^{-\alpha}
 \Big).
$$
Hence the relation (\ref{estmasse}) is equivalent to the following equality
$$ e^{\alpha}\Big([-2\alpha+1+e^{\alpha}+\delta(p)]e^{-\delta(p)}-e^{-2\alpha}-e^{-\alpha} \Big)+
e^{-\alpha}[\delta(p)-\ln(Y(p))-2\alpha]+Y(p)=1.$$
As $\delta(p)$ is chosen such that $f(\delta(p),Y(p))=0$, 
this conclude the proof of Theorem \ref{perio1} assuming that Lemma \ref{estf}holds. \hfill$\square$

\vspace{0,5cm}

\noindent {\bf Proof of Lemma \ref{estf}.}
We argue by continuity in $\delta$ and consider the endpoints. On the one hand, we have
$$
f(\alpha+p,Y(p))=\frac{e^{-\alpha}}{Y(p)}[1+e^{\alpha}+\ln(Y(p))]-e^{-\alpha}-e^{-2\alpha}+Y(p)-1.
$$
We are going to prove that 
$$
f(\alpha+p,Y(p)) \geq h(e^{-\alpha}) \geq 0.
$$
with the function $h:[e^{-\alpha},1] \to \RR$ given by 
$$
h(Y):= \frac{e^{-\alpha}}{Y}[1+e^{\alpha}+\ln(Y)]-e^{-\alpha}-e^{-2\alpha}+Y-1.
$$
To do so, we compute
$$
h'(Y)= 1+e^{-\alpha}\frac{Y- \ln(Y)}{Y^{2}}>0.
$$
We deduce the inequality.
\\

On the other hand,
$$
f(2\alpha,Y(p))=-e^{-\alpha}\ln(Y(p))+Y(p)-1 \leq 0,
$$
which concludes the proof of Lemma \ref{estf}. 
\hfill $\square$

\section{Simulations}
\label{sec:num}

This section deals with the numerical solution to the nonlinear renewal equation \eqref{as:id} with the transition from refractory state given by the piecewise smooth function \eqref{def_N} that has been used to derive analytical periodic solutions in section \ref{sec:periodic}. The two questions we address here are 1) which of these periodic solutions can we observe numerically, 2) what is the effect of smoothing terms that describe memory effects such as synaptic integration. In the present, synaptic integration occurs through a delay between the actual firing rate $N(t) $ defined as in \eqref{eqneuron} and its input in the firing rates coefficients $X(t)$ in section \ref{sec:model}. Rather than the simple rule \eqref{eq:interaction}, one introduces synaptic integration either through a convolution, as it was done in \cite{Pak_Per_Sal}, or through a differential equation which is easier for numerical simulations
\begin{equation}
\label{eq:delay}
\lb \frac{dX(t)}{dt} +X(t) = N(t).
\end{equation}
Here $\lb>0$ is interpreted as the relaxation time associated with the synaptic integration through a passive membrane. Our purpose is to see the effect of this parameter in the stability of the periodic solutions built in section \ref{sec:periodic}.
\\

\begin{figure}
{\centering
\includegraphics[width =7cm, height =6cm]{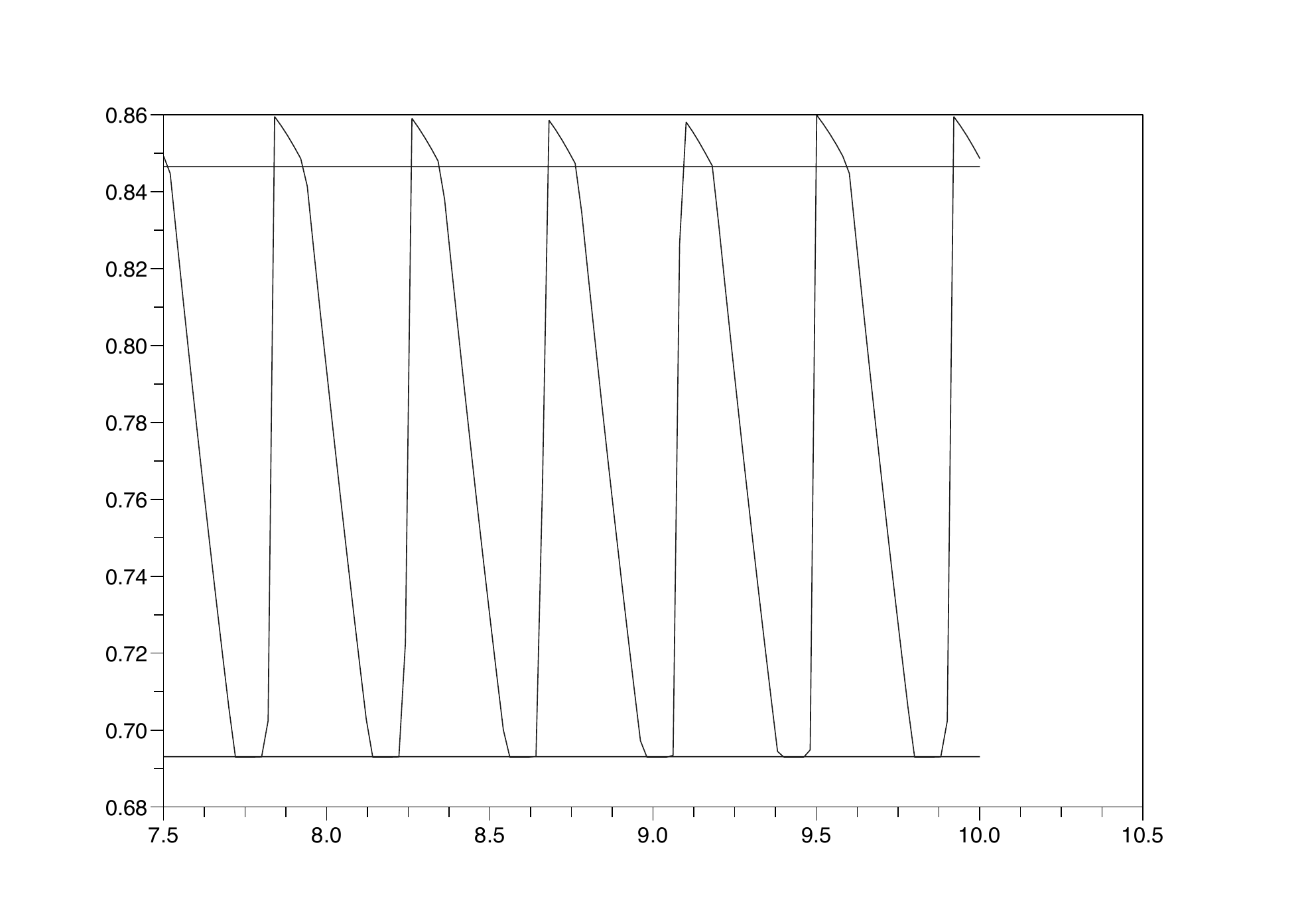} \hspace{.5cm}
\includegraphics[width =7cm, height =6cm]{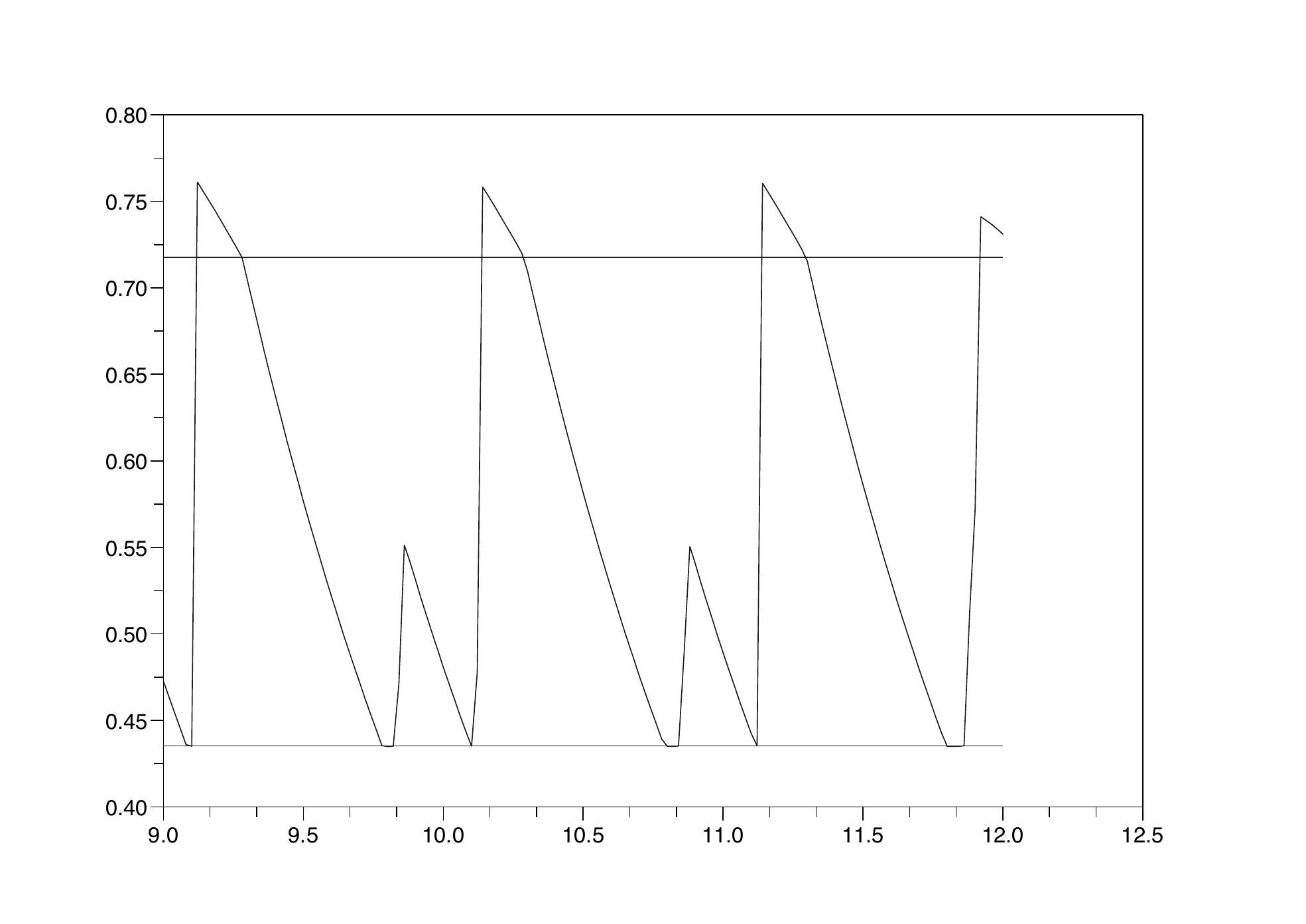}
\vspace{-.3cm} 
\caption{ \label{fig:nodelay1}  (Instantaneous transmission) Total neural activity $N(t)$ as computed with the numerical scheme \eqref{eq:numsch1}--\eqref{eq:numsch3} with the data in section \ref{sec:periodic}. Left: $\alpha=.2$. Right: $\alpha =.5$. The continuous lines give the values $N_-$ and $N_+$.  This is a numerical solution compatible with our construction of solutions in section \ref{sec:periodic3} with $p>0$.
} } \end{figure} 
\begin{figure}
{\centering
\includegraphics[width =7cm, height =6cm]{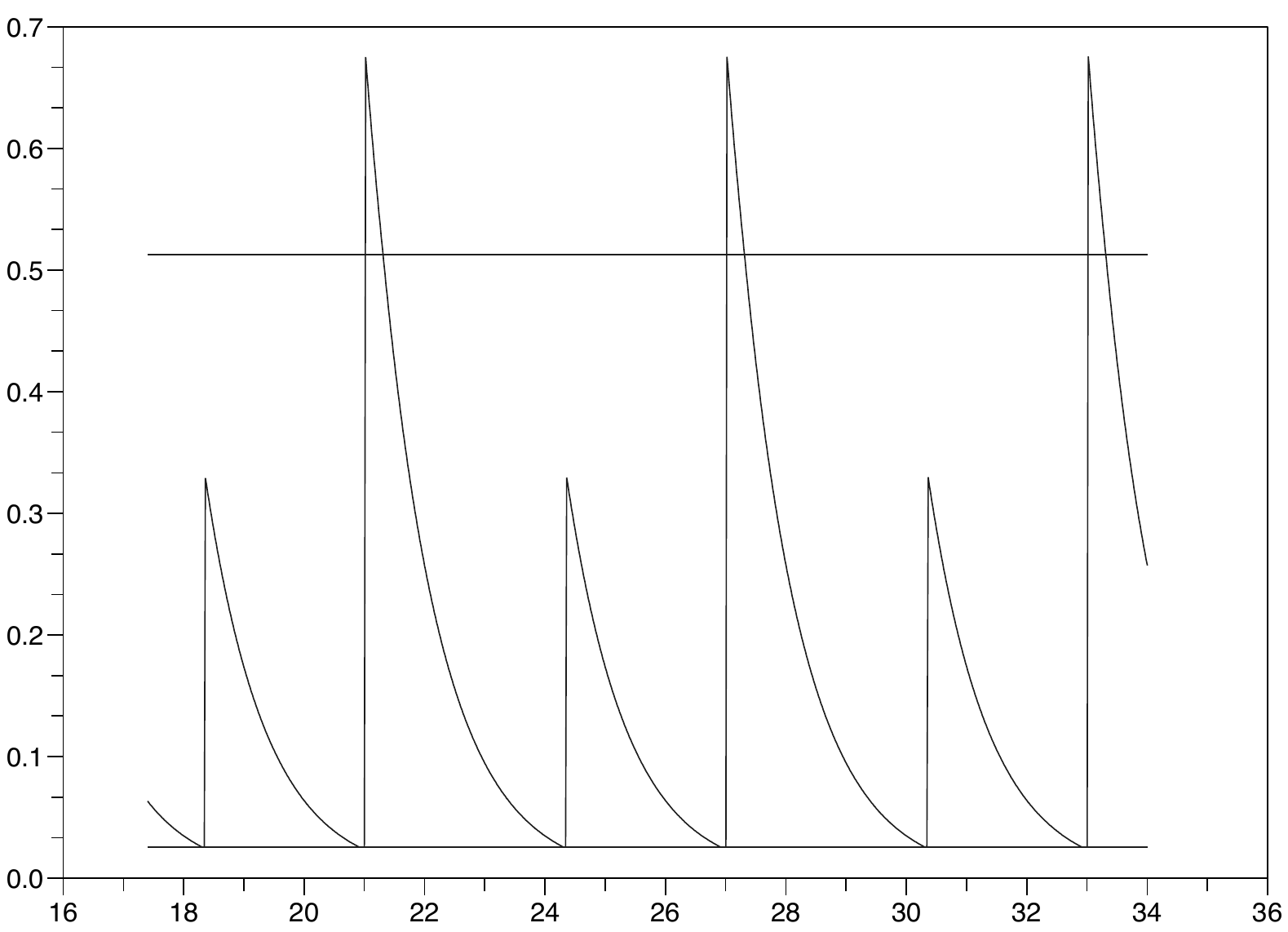} \hspace{.5cm}
\includegraphics[width =7cm, height =6cm]{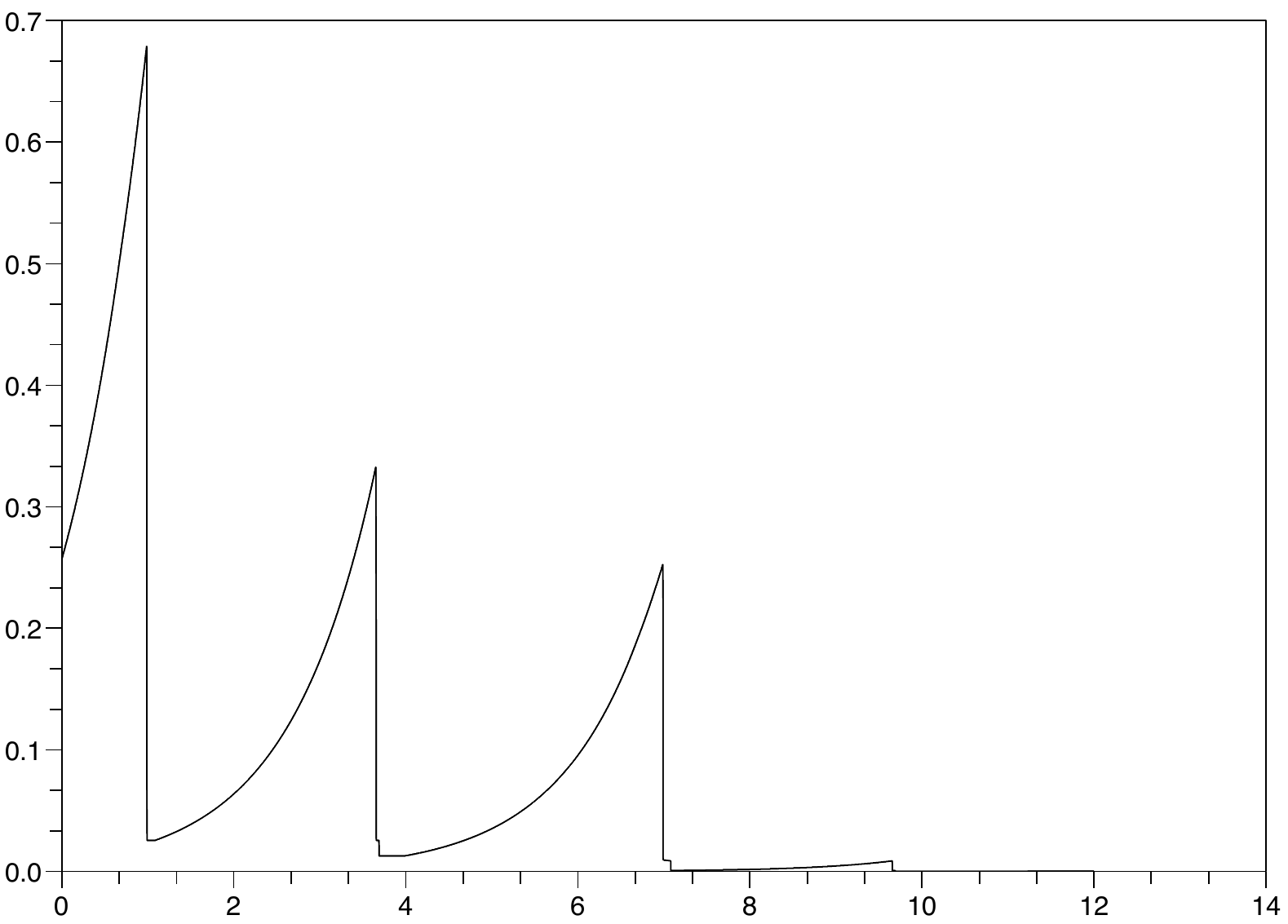}
\vspace{-.3cm} 
\caption{ \label{fig:nodelay2}  (Instantaneous transmission) Total neural activity $N(t)$ (left) and solution $n(s)$ (right)  computed as in Figure \ref{fig:nodelay1} with  $\alpha=3$. This numerical simulation is compatible with the exact solutions in section \ref{sec:periodic3} with $p>0$.
} } \end{figure} 
\begin{figure}
{\centering
\includegraphics[width =7cm, height =6cm]{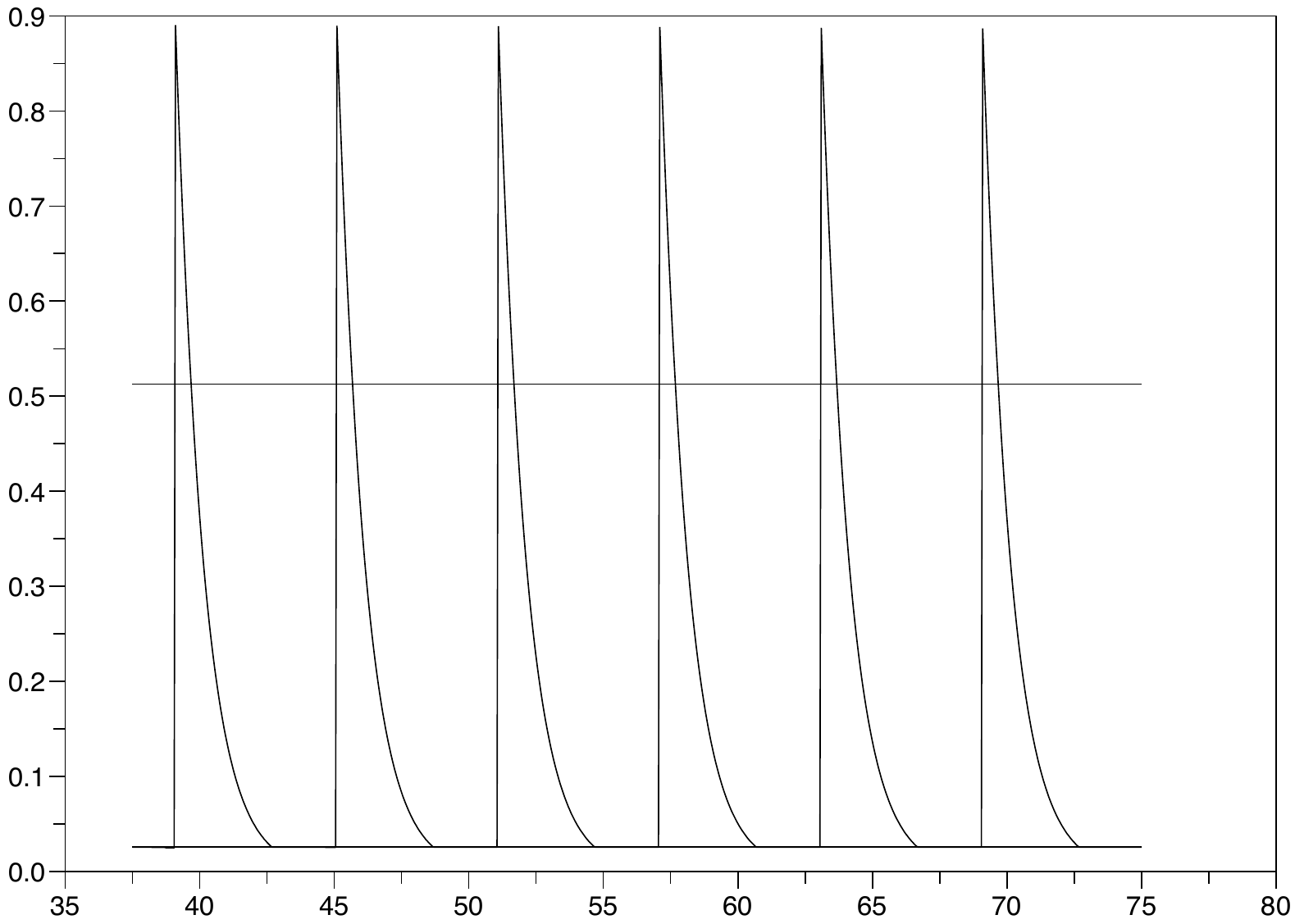} 
\vspace{-.3cm} 
\caption{ \label{fig:nodelay3}  (Instantaneous transmission) As in Figure \ref {fig:nodelay2} with the initial data $n^0(s)=e^{-s}$. Comparison shows that various periodic solutions are attractive depending on the initial data. Here we obtain periodic solutions as in section \ref{sec:periodic3} with $p=0$.
} } \end{figure} 

\subsection{Numerical methods}

We compute the numerical solution with a standard upwind scheme. The theory of such finite volume schemes can be found in \cite{bouchut,GR}  and in the references therein. Computationally it gives a very simple algorithm where the solution $n(s,t)$ is approximated by $n^k_i \approx \frac{1}{\Delta s} \int_{s_{i-1/2}}^{s_{i+1/2}} n(s,t^k) ds$ with grid parameters defined by $\Delta s=\Delta t$ (this means that the Courant-Friedrichs-Lewy (CFL in short) number equal to one, as it is usual for age-structured equations. We also use the notation 
$s_{i+1/2}=  (i+1/2) \Delta s $ and $t^k=k \Delta t$. The nonlinearities specific to the model \fer{eqneuron}, \fer{assimple} are treated as follows. The firing rate term $p\big(s,X(t)\big) n(s,t)$ is treated implicitly (and thus with  unconditional stability) which allows us to keep the maximal time step needed for accuracy in the transport term. Then one finds the discrete values $( n^k_i)_{1\leq i \leq I_M}$ iteratively on the time variable thanks to 
\beq \label{eq:numsch1}
\left\{\begin{array}{l}
\wt n^{k+1}_i = n^k_{i-1}/(1+\Delta t \; p_i^k),  \qquad 1\leq i \leq I_M,
\\
n^{k}_{0}:= N^{k} := \displaystyle\frac{1}{\Delta s} \sum_{i=1}^{I_M} p_i^k n^k_i ,
\end{array}\right.
\eeq
with $p_i^k= \ind{i \Delta s > \sg(X^k)}$. 
It has to be understood that the label $i=0$ stands for the boundary value ($s=0$ at the continuous level) and thus the formula \fer{eq:numsch1}. The semi-implicit scheme does not preserve the fundamental conservation law \fer{conservation} indicating we work with probabilities. Therefore one introduces a second order correction step (which in practice seems to have very little effect) and set 
\beq \label{eq:numsch2}
n^{k+1}_i =\wt n^{k+1}_i {\big /} \sum_{j=1}^{I_M} \wt n^{k+1}_i.
\eeq
When there is no delay we just use the idendity
\beq \label{eq:numsch3}
 X^{k+1} =N^{k}.
\eeq
When time delay is included according to \eqref{eq:delay}, we use an explicit Euler scheme and write 
\beq \label{eq:numsch4}
 X^{k+1} =X^{k} (1- \frac{\Delta t}{\lb})+ \frac{\Delta t}{\lb} N^{k}.
\eeq
Numerical solutions are computed with 1000 grid points per unit length (we have tested that this allows to reach numerical convergence). 

\subsection{Numerical periodic solutions}

The numerical results  are displayed in Figures \ref{fig:nodelay1} and \ref{fig:nodelay2} for times such that the periodic regime is established.  It seems the period is precisely   $2\alpha$. The numerical solution is therefore not the one peak periodic solution built in section  \ref{sec:periodic1} with period $\alpha$.  It is not either the periodic solution of period less than $2\alpha$ build in section  \ref{sec:periodic2}. This fact led us to build the third class in section  \ref{sec:periodic3} which precisely fits with the figures. The flat state can be observed even in figure \ref{fig:nodelay2} (zooming if necessary). 
\\

\begin{figure}
{\centering
\includegraphics[width =7cm, height =6cm]{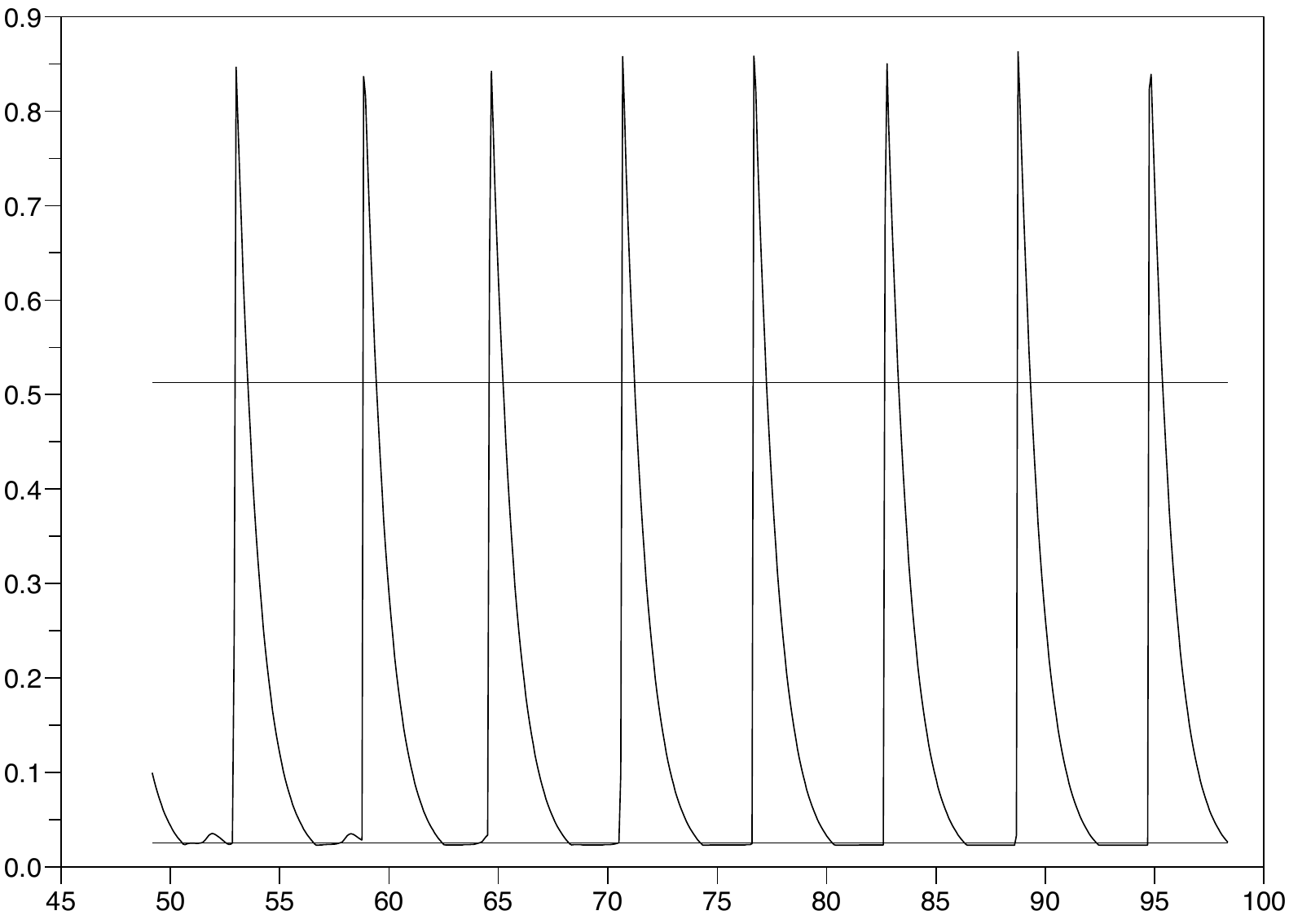} \hspace{.5cm}
\includegraphics[width =7cm, height =6cm]{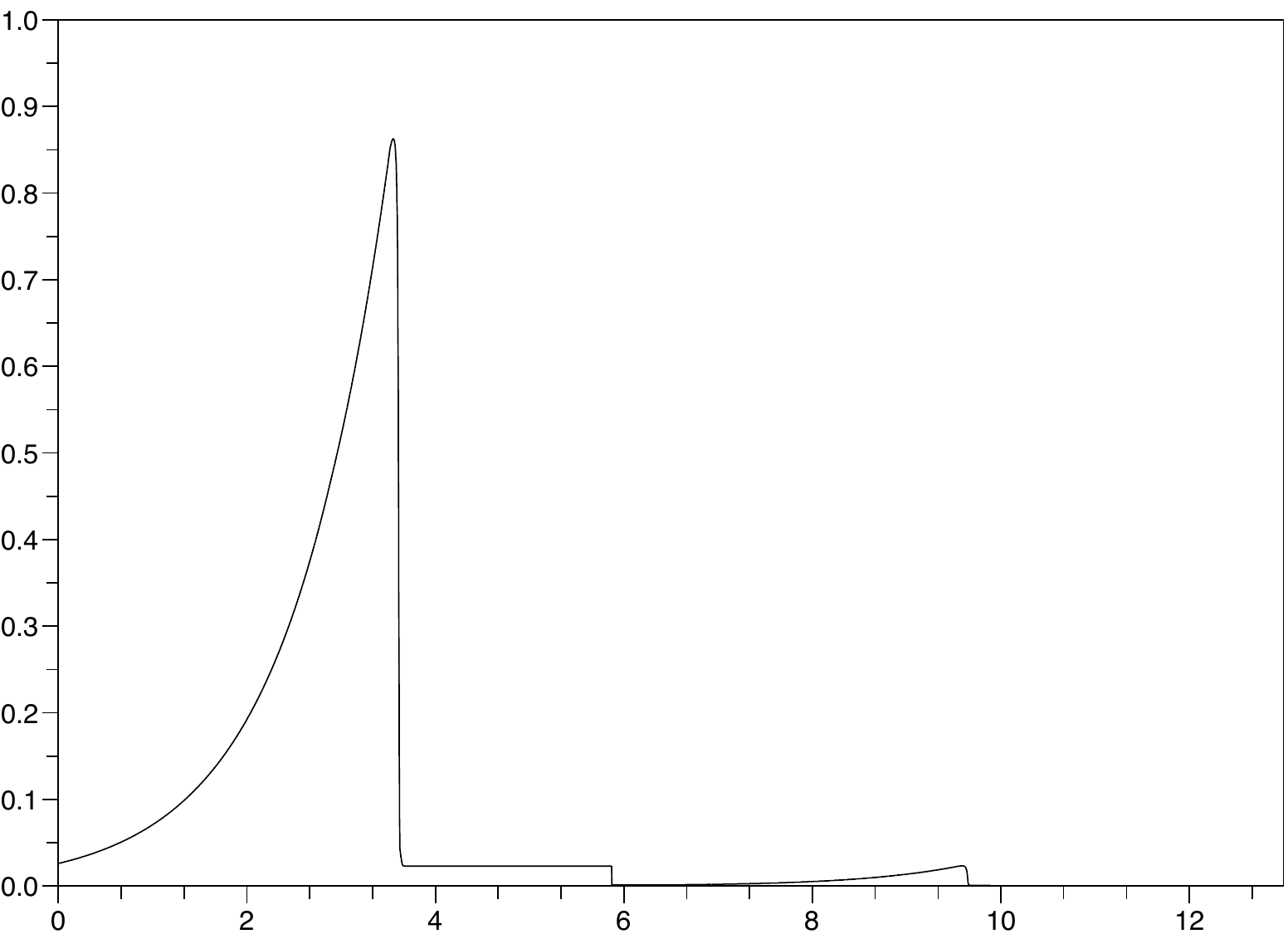}
\vspace{-.3cm} 
\caption{ \label{fig:delay}  (Time delay) Total neural activity $N(t)$ (left) and solution $n(s)$ (right)  computed as in Figure \ref{fig:nodelay1} with  $\alpha=3$. But time delay is included with $\lambda=.1$.
} } \end{figure} 

The effect of synaptic integration is to stabilize the dynamics and periodic solutions do not appear when $\lb$ is too large. When periodic solutions are still computed, the effect is to smooth out  the solution $n(s)$ away from the discontinuity which is still present. 

\subsection{Transition from desynchronization to sustained activity}

We have also included a parameter $J$ that describes the connectivity of the network. This can be implemented by simply using the threshold $\sigma(J N(t))$ in place of $\sigma(N(t))$ and 
$$
 p(s, N) = \ind{s> \sigma(J N))}.
$$
 
The numerical scheme shows a rapid transition at a critical value $J¬*$; for $J<J^*$ there is desynchronization and for $J>J^*$ discontinuous periodic solutions of large amplitudes appear similar to those depicted previously.  This confirms that the periodic solutions do not result from a supercritical Hopf-bifurcation.

\section{Conclusion}
\label{sec:conclusion}

The 'time elapsed' model for neural networks is a model that uses the probability density of neurons structured by the time elapsed since the last discharge and that has been introduced in \cite{Pak_Per_Sal}   and is based on the stochastic simulations in \cite{cf:pham98a, cf:pham98b}.
Following the study in \cite{Pak_Per_Sal}, both desynchronization or periodic solution may occur depending on the nonlinearity. This nonlinearity  takes into account the connectivity of the network through a modulation of the refractory period depending upon the total activity of the network.

Here we have quantified the nonlinearities that lead to total desynchronization. We arrive to a condition which is somehow sharp because we can build periodic solutions when the previous condition is not fulfilled. Numerical simulations show that the periodic solutions are numerous and depend on the initial data and the analytical solutions built in section \ref{sec:periodic} are not always the stable ones observed numerically.

We conjecture that the periodic solutions are always discontinuous (both the total network activity $N(\cdot)$ and the probability distribution of neurons $n$) in accordance with the numerical computations shown in section \ref{sec:num} and with the theoretical construction in section \ref{sec:periodic}.


\end{document}